\documentclass[10pt]{article}
\usepackage{amsmath,amsfonts,amsthm,amssymb,amscd}

\catcode`\é=\active\def é{\'e} \catcode`\è=\active\def è{\`e}
\catcode`\ç=\active\def ç{\c{c}} \catcode`\à=\active\def à{\`a}
\catcode`\ê=\active\def ê{\^e} \catcode`\ù=\active\def ù{\`u}
\catcode`\ô=\active\def ô{\^o} \catcode`\û=\active\def û{\^u}
\catcode`\î=\active\def î{\^{\i}} \catcode`\â=\active\def â{\^a}
\catcode`\ö=\active\def ö{\"o}

\newcommand{\be}{\begin{equation}}
\newcommand{\ee}{\end{equation}}
\newcommand{\bs}{\begin{split}}
\newcommand{\es}{\end{split}}
\newcommand{\ba}{\begin{align}}
\newcommand{\ea}{\end{align}}
\newcommand{\basl}[1]{\begin{align}\begin{split}\label{#1}}
\newcommand{\bas}{\begin{align}\begin{split}}

\newtheorem{theo}{Theorem}[section]
\newtheorem{prop}[theo]{Proposition}

\newtheorem{coro}[theo]{Corollary}
\newtheorem{rem}[theo]{Remark}
\newtheorem{exa}[theo]{Example}

\newcommand\N{\mathbb{N}}
\newcommand\Z{\mathbb{Z}}

\newcommand\R{\mathbb{R}}

\title{The Weyl symbol of Schrödinger semigroups }
\author{\sc
L. Amour, L. Jager and J. Nourrigat
}

\date{\it Université de Reims}

\begin{document}

\maketitle

\begin{abstract}
\noindent
In this paper, we study the Weyl symbol of the Schrödinger semigroup $e^{-tH}$, $H=-\Delta+V$, $t>0$, on $L^2(\R^n)$, with nonnegative potentials $V$ in $L^1_{\rm loc}$. Some general estimates like the $L^{\infty}$ norm concerning the symbol $u$ are derived. In the case of large dimension, typically for nearest neighbor or  mean field interaction potentials, we prove estimates with parameters independent of the dimension for the derivatives $\partial_x^\alpha\partial_\xi^\beta u$. In particular, this implies that the symbol of the Schrödinger semigroups belongs to the class of symbols introduced in \cite{AJN} in a high-dimensional setting. In addition, a commutator estimate concerning the semigroup is proved.
\end{abstract}

\parindent=0pt
\parskip 10pt
\baselineskip 15pt



\parindent = 0 cm

\section{Introduction.}
Let  $V$  be a nonnegative function in $L^2_{\rm loc}(\R^n)$. It is known that
$H= -\Delta + V(x)$ is essentially selfadjoint on
$C_0^{\infty}(\R^n)$ and we also denote by $H$ its unique selfadjoint extension. We may also suppose that $V$ is only in
$L^1_{\rm loc}(\R^n)$ and use Theorem X.32 in \cite{RS} to define $H$ as a selfadjoint operator with a suitable domain. In this paper, we are interested in the Weyl symbol $u(\cdot, t)$ of $e^{-tH}$, for each $t>0$. Since this operator is bounded
in $L^2(\R^n)$, its Weyl symbol is a priori a tempered distribution
 $U(t)$ on $\R^{2n}$ which satisfies,
\be\label{1}< e^{-tH} f , g> = < U(t) , H(f , g , \cdot) >, \ee
for all
$f$ and $g$ in ${\cal S}(\R^n)$, where $H(f , g, x , \xi)$ is the Wigner function  (c.f. \cite{CR} or \cite{L}),
\be\label{2}H(f , g, x , \xi) = (2\pi)^{-n} \int _{\R^{n}} e^{iv\cdot
\xi} f\left ( x - \frac{v}{2} \right ) \overline {g\left ( x +
\frac{v}{2}  \right )}  dv. \ee

The aim of this work is to study this Weyl symbol when $V$
is a $C^{\infty}$ potential describing a large number of particles in interaction, either for a nearest neighbor interaction model in a lattice, or for a mean field approximation model. Our hypotheses on the interaction potentials will in particular imply that $u(\cdot, t)$, the symbol of $e^{-tH}$,
is a $C^{\infty}$ function on $\R^{2n}$ and we shall give
estimates for its derivatives showing that $u(\cdot, t)$ belongs to a class
of symbols allowing a Weyl calculus in large dimension. This calculus is developed in \cite{AN} for the composition of symbols and in \cite{AJN} for norm estimates, where the constants implied in the inequalities are independent of the dimension.

Let us specify this class of symbols. In \cite{AN} and \cite{AJN}, we say that a continuous function $F$  on
$\R^{2n}$ is in $S_m (M, \rho, \delta)$ (where $m$ is an integer
$\geq 0$, $M\geq 0$, and $\rho$ and $\delta$ are two sequences
$(\rho_j)_{(j\leq n)}$ and $(\delta_j)_{(j\leq n)}$ of real numbers $\geq
0$) if, for all multi-indices $\alpha$ and $\beta$ in $\N^n$ satisfying $0\leq \alpha_j,\beta_j\leq m$, the derivative $\partial_x^{\alpha}
\partial_{\xi}^{\beta}F$ is a continuous and bounded function verifying,
\be\label{3} \Vert \partial_x^{\alpha} \partial_{\xi}^{\beta} F\Vert _{L^{\infty}
(\R^{2n})} \leq M \prod _{j\leq n} \rho_j^{\alpha_j}
\delta_j^{\beta_j}.\ee
Then, for functions belonging to the classes $S_m (M , \rho , \delta)$, we develop in \cite{AN} and \cite{AJN} a Weyl calculus in large dimension.

In the first section, the symbol $u(\cdot,t)$ of the semigroup is proved to be in $L^\infty(\R^{2n})$ with $|u(\cdot,t)|\leq 1$ almost
everywhere and additionally, $\int_{\R^n}u(x,\xi,t)dx$ and
$\partial_\xi^\beta u(\cdot,t)$ are estimated, without further regularity hypothesis on the interaction potential $V$. In the second section, we consider the semigroup in a large dimension setting with regular potentials and obtain estimates on all the derivatives of the symbol proving in particular that it lies in the class $S_m(M,\rho,\delta)$ above. Supplementary assumptions on potentials regarding the large dimension are naturally necessary at this step. Then, in the third section, two examples of Schrödinger semigroups in large dimension, satisfying the assumptions of section 2, are considered, namely, the nearest neighbor and the mean field approximation potentials. Moreover, a commutator property is also proved.

\section{First properties of the symbol of the semigroup.}

The first step consists in writing  the Weyl symbol of
$e^{-tH}$ with the Feynman Kac formula, under rather general hypotheses on the potential $V$. We make the choice to not first express the symbol with the distribution kernel, in order to avoid the use of Brownian bridges. Let $T>0$ and $n$ be an integer $\geq 1$. We denote by $B$ the Banach space
of continuous functions $\omega $  on $[0, T]$ taking values into $\R^n$ and vanishing at $t=0$. This space is endowed with the supremum norm, with the Borel $\sigma-$algebra ${\cal B}$ and with the  Wiener measure
$\mu$ of variance 1 (c.f. Kuo \cite{K}).

\begin{prop}\label{p1} \it Let $V\geq 0$ be a function in
$L^1_{\rm loc}(\R^n)$. Let $U( t) $ be the Weyl symbol of the operator
$e^{-tH}$, first considered as a tempered  distribution on
$\R^{2n}$. Then, $U(t)$ is identified with a function $u(\cdot , t)$
in $L^{\infty}(\R^{2n})$. We have, for  each $t$ in $(0, T]$ and
for almost every $(x , \xi)$ in $\R^{2n}$,

\be\label{4} u(x , \xi, t)   = \int _{B }  e^{-i \omega (t) \xi} e^{-\int
_0^t V (x -\frac{\omega (t)}{2} +  \omega (s)) ds }d\mu (\omega ).
\ee
Moreover, the following inequality holds,
\be\label{5} |u(\cdot, t)|\leq 1, \ee
almost everywhere on $\R^{2n}$, for each $t\in[0,T]$.
\end{prop}

One notices that the above integral involves all of the trajectories of the Brownian motion with a starting and a finishing point that are symmetric with respect to $x$.

{\it Proof.} Let $f$ and $g$ in ${\cal S}(\R^n)$. When $V\geq 0$
belongs to $L^1_{\rm loc}(\R^n)$, one may apply
Feynman Kac formula (c.f. B. Simon \cite{S1} or \cite{S2}) written as,
\be\label{6}< e^{-tH} f , g> = \int _{\R^n\times B}  f( x + \omega (t))
\overline {g(x)}  e^{-\int _0^t V (x + \omega (s)) ds } dxd\mu
(\omega ).\ee
According to the Wigner function definition, we have for all
$x$ and $y$ in $\R^n$,
$$f(x) \overline {g(y)} = \int _{\R^n} H \left ( f, g , \frac{x+y}{2}
, \xi \right ) e^{-i (x-y) \cdot \xi} d\xi.$$
Consequently, for all $\omega $ in $B$,
$$f(x+ \omega (t) ) \overline {g(x )} =\int _{\R^n} H \left ( f, g ,
x + \frac{\omega (t)}{2} , \xi \right ) e^{-i\omega (t) \cdot \xi }
d\xi.$$
The Weyl symbol $U(t)$ of $e^{-tH}$ being a priori defined as a tempered distribution  on $\R^{2n}$, thus satisfies, for all
$F$ in ${\cal S}(\R^{2n})$,
$$ < U(t), F> = \int _{\R^{2n}\times B} F \left (
x + \frac{\omega (t)}{2} , \xi \right )e^{-i\omega (t) \cdot \xi }
e^{-\int _0^t V (x + \omega (s)) ds }  dx d\xi d\mu (\omega).$$
The above identity shows that, for all $F$ in ${\cal
S}(\R^{2n})$,
$$ |< U(t), F>| \leq \Vert F\Vert _{L^1(\R^{2n})}.$$
As a consequence, $U(t)$ is identified with a function $u(\cdot, t)$ in
$L^{\infty }(\R^{2n})$, with a $L^{\infty }(\R^{2n})$ norm
smaller or equal than $1$, and satisfying (\ref{4}) and (\ref{5}). The proposition is then
proved.\hfill$\Box$

As a first consequence of Proposition \ref{p1}, we give below two corollaries which do not assume that the potential $V$ is differentiable.

\begin{coro}\it For every multi-index $\beta$, for each
$t\geq 0$, the derivative $\partial _{\xi}^{\beta}u(x , \xi, t)$ understood in the sense of distributions, is a function in
$L^{\infty}(\R^{2n})$, which satisfies,
\be\label{7}\Vert \partial _{\xi}^{\beta} u(\cdot , t)\Vert _{L^{\infty}(\R^{2n})} \leq
 t^{|\beta|/2} \prod _{j\leq n} A_{\beta_j}, \hskip 2cm A_k =\frac{2^{k/2}}{ \sqrt {\pi}} \Gamma((k+1)/2).
\ee
Let $m$ be an integer $\geq 0$. If the multi-index $\beta$ verifies
$\beta_j \leq m$ for all $j\leq n$, then we have,
\be\label{8}\Vert \partial _{\xi}^{\beta} u(\cdot , t)\Vert _{L^{\infty}(\R^{2n})}
\leq B_m^{|\beta|} t^{|\beta|/2}, \hskip 2cm B_m = \max _{k\leq m}
A_k. \ee
When $m\geq 4$, we have $B_m = A_m$.
\end{coro}

{\it Proof.} We use the notation $\omega (t) = ( \omega_1 (t) , ... , \omega_n
(t))$. In view of Proposition \ref{p1},
$$ \Vert \partial _{\xi}^{\beta} u(\cdot , t)\Vert _{L^{\infty}(\R^{2n})} \leq
\int _B \prod _{j\leq n} |\omega_j (t)|^{\beta_j} d\mu (\omega).
$$
According to Kuo \cite{K}, we know that,
\be\label{9}\int _B \prod _{j\leq n} |\omega_j (t)|^{\beta_j} d\mu (\omega) =
 t^{|\beta|/2} \prod _{j\leq n} A_{\beta_j}, \ee
where the $A_k$ are given in (\ref{7}). This proves the corollary.\hfill$\Box$

\begin{coro}  If $V\geq 0$  and if the right hand side below defines a convergent integral, then, for each $t>0$ and for almost every  $\xi$ in $\R^n$, the function
$u(\cdot , \xi, t)$ belongs to $L^1(\R^n)$, and we have,
\be\label{10}\int _{\R^n} |u(x , \xi, t)|  dx \leq \int _{ \R^n}
e^{-tV(x)} dx. \ee
\end{coro}

{\it Proof.} From $(\ref{4})$, we see
$$|u(x , \xi, t)| \leq \int _{B }  e^{-\int
_0^t V (x -\frac{\omega (t)}{2} +  \omega (s)) ds }d\mu (\omega ).
$$
Since the function $x\mapsto e^{-tx}$ is convex, using Jensen inequality and integrating over $x\in\R^n$, for almost every $\xi\in\R^n$, lead to inequality $(\ref{10})$.
The corollary is thus proved.\hfill$\Box$

\section{The large dimension setting.}
We shall here give
 Hamiltonians
$H_{\Lambda}$ for systems with a large number of particles indexed on  $\Lambda$, for which we shall obtain estimates on the
derivatives of the Weyl symbol $u_{\Lambda}(\cdot, t)$ of $e^{-t H_{\Lambda}}$. These estimates prove in particular that
$u_{\Lambda}(\cdot , t)$ belongs to the class of symbols studied
in \cite{AN} and \cite{AJN},  allowing a Weyl calculus where all the
constants in the inequalities are independent of $\Lambda$. The assumptions on the interaction potentials $V_\Lambda$ are stated below. In the next section, we shall give two examples of Hamiltonians satisfying these hypotheses.

We suppose that the functions
$V_{\Lambda}$ are given, $\geq 0$, $C^{\infty}$  on $\R^{\Lambda}$, for each finite
subset  $\Lambda $ in $\Gamma$, for a given infinite countable set $\Gamma$. For  all
integers $m\geq 0$, we denote by ${\cal M}_m (\Lambda)$ the set of
multi-indices $\alpha$ in $\N^{\Lambda}$ such that $0\leq \alpha_j
\leq m$, for all $j\in \Lambda$. For each multi-index $\alpha$,
$S(\alpha)$ denotes the set of sites $j\in \Lambda$ such that
$\alpha_j \not= 0$.

Set $m\geq 1$. We also assume that there exists $C_m>0$
such that, for all finite subsets  $\Lambda $ in $\Gamma$,  for
all $\alpha $ in ${\cal M}_m (\Lambda)$, we have for all $x\in\R^\Lambda$
\be\label{11} \sum _{0\not =\beta \leq \alpha}|\partial ^{\beta } V_{\Lambda } (x)|
\leq C_{m}  |S(\alpha)|.\ee
We set,
\be\label{12}H_{\Lambda} = - \Delta_{\Lambda } + V_{\Lambda}(x) \ee
and   $U_{\Lambda}(t)$ denotes the Weyl symbol of
$e^{-tH_{\Lambda}}$ which, according to Proposition \ref{p1}, is a tempered distribution on
$\R^{\Lambda}\times \R^{\Lambda}$ identified with a function
$u_{\Lambda} (\cdot , t)$ in $L^{\infty}(\R^{\Lambda}\times
\R^{\Lambda})$.

\begin{theo}\label{t4} With these notations, let the functions
$V_{\Lambda}\geq 0$ in $C^{\infty }(\R^{\Lambda})$ be given for
all finite subsets $\Lambda$ of $\Gamma$. Let $m\geq 1$. We suppose
that there exists $C_m>0$ independent of
$\Lambda$, such that (\ref{11}) is satisfied.  For each $t>0$, and for
every finite subset $\Lambda$ of $\Gamma$, let $u_{\Lambda}
(\cdot , t)$ be the Weyl symbol of $e^{-tH_{\Lambda}}$, which is identified with a function in $L^{\infty} (\R^{\Lambda}\times
\R^{\Lambda})$ in view of Proposition \ref{p1}. Then, for each $\alpha$
and $\beta$ in ${\cal M}_m (\Lambda)$, the derivative
$\partial_x^{\alpha}
\partial_{\xi}^{\beta} u(\cdot , t)$, understood in the sense of
distributions, is a function in $L^{\infty}(\R^{\Lambda}\times
\R^{\Lambda})$ which satisfies,
\be\label{13}\Vert \partial_x^{\alpha} \partial_{\xi}^{\beta} u(\cdot , t)\Vert _{L^{\infty}
(\R^{\Lambda}\times \R^{\Lambda})} \leq m ^{|S(\alpha)|} e^{ t C_m
|S(\alpha)|}B_m ^{|S(\beta)|} t^{|\beta|/2}, \ee
where  $B_m$ is defined in (\ref{7}) and (\ref{8}), and $C_m$ in (\ref{11}) (these constants are independent of $\Lambda$).
\end{theo}

{\it Proof.} According to Proposition \ref{p1}, we have for all
$F$ in ${\cal S}(\R^{\Lambda}\times \R^{\Lambda})$,
\be\label{14}|< \partial_x^{\alpha }\partial_{\xi}^{\beta} U_{\Lambda}(t) , F>|\leq \Vert F \Vert
_{L^1(\R^{\Lambda}\times \R^{\Lambda})}\sup _{(x, \omega)\in \R^{\Lambda} \times B}  \left |
\partial_x^{\alpha }e^{-\int
_0^t V_{\Lambda}  (x + \omega (s)) ds } \right |\int _B \prod _{j\in
\Lambda} |\omega_j (t)|^{\beta_j} d\mu (\omega).
 \ee
We shall use a multi-dimensional variant of Faà di
Bruno formula due to Constantine Savits \cite{CS}. For each multi-index
$\alpha$, denote by $F(\alpha)$ the set of mappings $\varphi$
from the set of multi-indices $0\not =\beta \leq \alpha$ into the set of
integers $\geq 0$, such that
$$ \sum _{0\not =\beta\leq \alpha } \varphi (\beta) \beta = \alpha. $$
 Constantine Savits formula is rewritten as,
\be\label{15}\partial^{\alpha}e^{W (x)} = \alpha ! e^{W
(x)}  \sum _{\varphi \in F(\alpha)} \prod _{0\not =\beta \leq \alpha}
 \frac{1
}{\varphi (\beta)!} \left [ \frac{\partial ^{\beta }W(x)}{\beta
! }\right ]^{\varphi (\beta) }. \ee
For each $t>0$ and for almost all $\omega$ in $B$, we apply this formula with
$$ W(x) =- \int
_0^t V_{\Lambda}  (x + \omega (s)) ds. $$
Since $V_{\Lambda}\geq 0$, we obtain
$$\sup _{(x, \omega)\in \R^{\Lambda} \times B}  \left |
\partial_x^{\alpha }e^{-\int
_0^t V_{\Lambda}  (x + \omega (s)) ds } \right |  \leq \alpha ! \sum
_{\varphi \in F(\alpha)} \prod _{0\not =\beta \leq \alpha} \frac{1 }{\varphi
(\beta)!}  \left [ \frac{t \Vert \partial ^{\beta }V_{\Lambda}\Vert
_{L^{\infty} }}{ \beta ! } \right ]^{\varphi (\beta) }.  $$
Besides,
$$\sum
_{\varphi \in F(\alpha)} \prod _{0\not =\beta \leq \alpha} \frac{1  }{\varphi
(\beta)!}  \left [ \frac{t \Vert \partial ^{\beta }V_{\Lambda}\Vert
_{L^{\infty} } }{ \beta ! } \right ]^{\varphi (\beta) } \leq \exp
\left [ \sum _{0\not =\beta \leq \alpha } \frac{t \Vert\partial ^{\beta
}V_{\Lambda}\Vert _{L^{\infty} } }{ \beta!} \right ].$$
The last factor in (\ref{14}) is bounded using (\ref{9})(\ref{8}) and the above right hand side is bounded using the hypothesis (\ref{11}). We then deduce that,
$$ |< \partial_x^{\alpha }\partial_{\xi}^{\beta}U_{\Lambda}(t) , F>|\leq \alpha ! \Vert F \Vert
_{L^1(\R^{2n})} e^{ t C_m |S(\alpha)|} B_m ^{|S(\beta)|}
t^{|\beta|/2}.  $$
Since $\alpha$ is in ${\cal M}_m (\Lambda)$ , we have $\alpha !
\leq m ^{|S(\alpha)|}$. The proof of Theorem \ref{t4} is then completed.
\hfill$\Box$

\begin{rem}
Theorem \ref{t4}
shows that, if the family of functions $(V_{\Lambda})$ verifies (\ref{11}) with $C_m>0$, then, for all $t>0$, and for each $m\geq 0$,
the family of functions $u_{\Lambda}(\cdot , t)$ belongs to the class
$S_m (1 , \rho , \delta)$ defined in (\ref{3}), where $ \rho_j = m
e^{tC_m}$ and $\delta_j = B_m \sqrt {t}$ for all $j\in \Gamma$.
We remark that $\rho_j$ and $\delta_j$ depend on $m$ but not on $\Lambda$, and thus, not on the dimension.
\end{rem}

\section{Examples and application.}
\subsection{Two examples.}

We shall in  this section give two examples of family of potentials
 $(V_{\Lambda})$ satisfying (\ref{11}) for all integers $m\geq 1$. The first one corresponds to the nearest neighbor interaction in a lattice and the second one corresponds to the mean field approximation model.

\begin{exa}\label{e1} Set $\Gamma =\Z^d$ ($d\geq 1$). Let $F$ and $G$ be two nonnegative functions in $C^{\infty}(\R)$,
bounded together with all their derivatives. For each
finite subset $\Lambda$ of $\Gamma$, we set,
\be\label{16}V_{\Lambda}(x) = \sum _{j\in \Lambda} F(x_j) + \sum _{(j , k)\in
\Lambda^2 \atop |j-k|_{\infty}  = 1} G( x_j - x_k). \ee
Then, for all integers $m\geq 1$, there exists $C_m>0$ such that the
family of  functions $(V_{\Lambda})$ satisfies (\ref{11}).
\end{exa}

\begin{exa}\label{e2}  Let $\Gamma $ be an infinite countable set. Fix a function $G\geq 0$ in
$C^{\infty}(\R)$, bounded together with all its derivatives.  Let, for
each finite subset $\Lambda$ of $\Gamma$,
\be\label{17} V_{\Lambda}(x) =  \frac{1}{ |\Lambda|} \sum _{(j , k)\in
\Lambda^2} G( x_j - x_k). \ee
Then, for any integer $m\geq 1$, there is $C_m>0$ such that the
family of potentials $(V_{\Lambda})$ is verifying (\ref{11}).
\end{exa}

\subsection{Application.}

For all finite subsets $\Lambda$ in $\Gamma$, choose a
function  $p_{\Lambda}\geq 0$  in the Schwarz space $\mathcal{S}(\R^{\Lambda} \times
\R^{\Lambda})$. It is known that the Weyl operator
${\rm Op}^{\rm Weyl}(p_{\Lambda})$ is trace class. We  suppose that its trace
equals 1.

For all finite subsets $\Lambda$ in $\Gamma$, suppose that the
functions $V_{\Lambda}\geq 0$ in $\R^{\Lambda}$ are given satisfying the conditions of
Example \ref{e1} or Example \ref{e2}, and denote by $H_{\Lambda}$ the Hamiltonian defined in (\ref{12}). Let $A$ be a function on $\R$ and denote by $A_j $  the  multiplication operator by
the function $A(x_j)$ ($j\in \Lambda$). The function $A$ is chosen to be polynomial to avoid a long development on pseudodifferential operators.

\begin{prop} With these notations, there exists a constant
$C>0$ such that, for all finite subsets $\Lambda$ of $\Gamma$,
for each $t$ in $(0, 1]$, for every $j$ in $\Lambda$, we have,
$$ \Big |{\rm Tr} ( [ A_j , e^{-tH_{\Lambda}} ] \circ
{\rm Op}^{\rm Weyl}(p_{\Lambda}))
\Big | \leq C \sqrt {t}.$$
\end{prop}

{\it Proof.} Let $F_{\Lambda , t}$ be the Weyl symbol of the
commutator $[ A_j , e^{-tH_{\Lambda}} ]$. According to  Theorem
\ref{t4} and to the Weyl calculus in one dimension, there exists $C>0$ such that, for
all finite $\Lambda$ in $\Gamma$, and for any $t$ in $(0, 1]$,
$$\Vert F_{\Lambda , t}\Vert _{L^{\infty}
(\R^{\Lambda}\times \R^{\Lambda})} \leq  C \sqrt {t}.$$
It is known that,
$${\rm Tr} ( [ A_j , e^{-tH_{\Lambda}} ] \circ
{\rm Op}^{\rm Weyl}(p_{\Lambda}))  = (2\pi)^{-|\Lambda|} \int
_{\R^{\Lambda}\times \R^{\Lambda}} F_{\Lambda , t}(x , \xi)
p_{\Lambda} (x , \xi) dx d\xi, $$
$${\rm Tr} ({\rm Op}^{\rm Weyl}(p_{\Lambda})) =(2\pi)^{-|\Lambda|} \int
_{\R^{\Lambda}\times \R^{\Lambda}}p_{\Lambda} (x , \xi) dx d\xi =1.
$$
The proposition then follows.\hfill$\Box$

\medskip
{\it E-mail address: }\{laurent.amour, lisette.jager, jean.nourrigat\}@univ-reims.fr

{\it Address:} LMR
EA 4535 and FR CNRS 3399, Universit\'e de Reims Champagne-Ardenne,
 Moulin de la Housse, BP 1039,
 51687 REIMS Cedex 2, France.

\end{document}